\documentclass[pdflatex,sn-mathphys-ay]{sn-jnl}

\usepackage{graphicx}
\usepackage{multirow}
\usepackage{amsmath,amssymb,amsfonts}
\usepackage{amsthm}
\usepackage{mathrsfs}
\usepackage[title]{appendix}
\usepackage{xcolor}
\usepackage{textcomp}
\usepackage{manyfoot}
\usepackage{booktabs}
\usepackage{indentfirst}
\usepackage{xr}

\theoremstyle{thmstyleone}
\newtheorem{theorem}{Theorem}

\newtheorem{coro}{Corollary}
\newtheorem{prop}{Proposition}

\theoremstyle{thmstyletwo}
\newtheorem{remark}{Remark}

\theoremstyle{thmstylethree}
\newtheorem{defi}{Definition}
\newtheorem{example}{Example}

\newcommand{\ba}{\begin{array}}
\newcommand{\ea}{\end{array}}
\newcommand{\bt}{\begin{tabular}}
\newcommand{\et}{\end{tabular}}
\newcommand{\btb}{\begin{table}}
\newcommand{\etb}{\end{table}}
\newcommand{\bc}{\begin{center}}
\newcommand{\ec}{\end{center}}
\newcommand{\bea}{\begin{eqnarray}}
\newcommand{\eea}{\end{eqnarray}}
\newcommand{\Bea}{\begin{eqnarray*}}
\newcommand{\Eea}{\end{eqnarray*}}
\newcommand{\beq}{\begin{equation}}
\newcommand{\eeq}{\end{equation}}

\newcommand{\SD}{{\rm SD}}

\newcommand{\rd}{{\rm d}}
\newcommand{\bu}{{\boldsymbol u}}
\newcommand{\bx}{{\boldsymbol x}}
\newcommand{\by}{{\boldsymbol y}}
\newcommand{\bz}{{\boldsymbol z}}
\newcommand{\bo}{{\boldsymbol \omega}}
\newcommand{\rSOA}{{\mathrm{SOA}}}
\newcommand{\rOA}{{\mathrm{OA}}}
\newcommand{\dist}{d}

\makeatletter
\let\saved@bibcite\bibcite
\let\bibcite\@gobbletwo
\let\bibcite\saved@bibcite
\makeatother

\begin{document}

\title[Uniform projection designs under the stratified $L_2$-discrepancy]{Uniform projection designs under the stratified $L_2$-discrepancy}

\author[1]{\fnm{Sixu} \sur{Liu}}\email{liusixu@bimsa.cn}
\author*[2]{\fnm{Yaping} \sur{Wang}}\email{ypwang@fem.ecnu.edu.cn}

\affil[1]{\orgname{Beijing Institute of Mathematical Sciences and Applications}, \orgaddress{\city{Beijing}, \postcode{101408}, \country{China}}}
\affil*[2]{\orgname{KLATASDS--MOE, School of Statistics, East China Normal University}, \orgaddress{\city{Shanghai}, \postcode{200062}, \country{China}}}

\abstract{This paper studies a uniform projection criterion for space-filling designs under the stratified $L_2$-discrepancy. The criterion, denoted by $\Phi_{\rm SD}$, is the average squared stratified $L_2$-discrepancy over all two-dimensional projections. For U-type $(n,m,s^p)$ designs, we derive an explicit formula for $\Phi_{\rm SD}$ in terms of row-pairwise weighted hierarchical distances, and we establish sharp lower and upper bounds with equality conditions. We further show that many known optimal constructions attain the lower bound of $\Phi_{\rm SD}$, and that designs attaining the lower bound of the full stratified $L_2$-discrepancy also attain the lower bound of $\Phi_{\rm SD}$. The criterion can be evaluated in $O(n^2m)$ time, with a modest reduction in arithmetic operations compared with direct projection-wise evaluation. Numerical studies illustrate the theoretical results and show that $\Phi_{\rm SD}$ is effective for assessing low-dimensional projection uniformity. }


\maketitle

\section{Introduction}
Space-filling designs have attracted growing attention for the investigation of complex systems in both computer and physical experiments~\citep{Fang06, Fang18, Santner03}.
 These designs aim to distribute points uniformly across the experimental domain, with discrepancy serving as a common criterion for evaluating their uniformity. In general, discrepancy quantifies how much the empirical distribution of a design deviates from the theoretical uniform distribution over the design space.

Among various measures of discrepancy, generalized $L_2$-type discrepancies have been widely adopted \citep{Hickernell98a, Hickernell98b, Hickernell02, ZhouFangNing13}. Introduced by \cite{Hickernell98a, Hickernell98b}, this family of discrepancies is constructed within a reproducing kernel Hilbert space and accounts for the projection properties of design points across all subdimensions. Notable members include the (modified) $L_2$-star discrepancy, the centered $L_2$-discrepancy (CD), the wrap-around $L_2$-discrepancy (WD), and the mixture $L_2$-discrepancy (MD).

A critical issue with commonly used discrepancies like the CD, as highlighted by \citet{ZhouFangNing13}, is their inherent dimensionality effect. In high-dimensional spaces, this phenomenon often leads to an over-concentration of design points toward the center of the experimental domain. Consequently, a design may achieve satisfactory global uniformity while exhibiting poor space-filling properties in its lower-dimensional projections. In many practical applications, system responses are often dominated by main effects and low-order interactions, a scenario known as effect sparsity; thus, ensuring the uniformity of low-dimensional projections is essential.

To mitigate these issues, \citet{SunWangXu19} introduced the \textit{uniform projection criterion}, which seeks to minimize the average CD across all possible two-dimensional projections. This approach builds upon foundational projection-based methodologies; for instance, \citet{Hickernell02} developed the projection discrepancy pattern for assessing uniformity in projection subspaces, while \citet{FangQin05} established the uniform pattern and minimum projection uniformity specifically for two-level designs.

The literature on {uniform projection designs} (UPDs) has expanded considerably in recent years. Regarding construction methodologies, UPDs have been developed using good lattice points \citep{SunWangXu19,WangSunXu22}, level permutation and expansion techniques \citep{Zhou23}, and tuned differential evolution algorithms \citep{onyambu2025tuning}. Parallel research has explored maximin $L_1$-distance U-type constructions \citep{WangXiaoXu18, LiLiuTang21}. Furthermore, the framework of projection uniformity has been generalized to encompass other $L_2$-type discrepancies \citep{Liu23}, adapted for nested Latin hypercube designs \citep{chen2021uniform}, and theoretically linked with strong orthogonal arrays \citep{sun2023uniform}.

 Recently, \citet{TianXu25} introduced the stratified $L_2$-discrepancy as a robust measure for evaluating design uniformity when the experimental domain is partitioned into subregions. Formulated within the reproducing kernel framework, this discrepancy inherits the mathematical advantages of generalized $L_2$-discrepancies while mitigating the concentration of design points in specific sub-domains---a critical issue in high-dimensional spaces.
 Designs that minimize this discrepancy are not only space-filling but also exhibit favorable projection properties.
 Furthermore, the stratified $L_2$-discrepancy is highly flexible, as it accommodates various stratification schemes and importance weights, enabling a more nuanced assessment of projection quality.

In this paper, we systematically investigate a uniform projection criterion that minimizes the average stratified $L_2$-discrepancy across all two-dimensional projections of a design. Compared to existing projection criteria based on generalized $L_2$-type discrepancies \citep{SunWangXu19, Liu23}, our proposed criterion provides enhanced stratification properties. We derive an analytical formula for the criterion through row-pairwise relationships and establish its explicit connection to a weighted hierarchical discrete distance. This theoretical framework allows us to provide sharp lower and upper bounds, along with the corresponding equality conditions.
  We further demonstrate that several established constructions, including those based on general strong orthogonal arrays (GSOAs) derived from multiplication tables and generalized Hadamard matrices \citep{TianXu24, TianXu25, Gao25}, achieve the theoretical lower bound. Notably, it is found that any design attaining the lower bound of the full stratified $L_2$-discrepancy \citep{TianXu25} also achieves the lower bound of our criterion, confirming that our approach successfully identifies designs that are optimal under the broader stratified framework. From a computational perspective, evaluating our criterion for a design with $n$ runs and $m$ factors requires $O(n^2m)$ operations. This matches the complexity of the full stratified $L_2$-discrepancy. Numerical studies further validate our theoretical findings and illustrate the effectiveness of the proposed criterion in assessing low-dimensional projection uniformity.

The remainder of the paper is organized as follows. Section~\ref{sec2} introduces notation and reviews the stratified $L_2$-discrepancy. Section~\ref{sec3} formally introduces the projection criterion $\Phi_{\rm SD}$, derives its fast-computing formula based on row-pairwise distances, and provides the associated lower and upper bounds.
Section~\ref{sec3opt} presents various constructions for optimal designs that achieve the lower bound of $\Phi_{\rm SD}$.
 Section~\ref{sec4} reports numerical studies to validate the performance of the proposed criterion.
  Section~\ref{sec5} explores further connections between our criterion, the space-filling hierarchy, and projection pattern enumerators. Section~\ref{sec:conclusion} concludes the paper with a brief discussion.

\section{Preliminaries}\label{sec2}

\subsection{Notation}
For a positive integer $q$, let $\mathbb{Z}_q=\{0,1,\ldots,q-1\}$.
The notation $(n, m, s)$ design refers to a design with $n$ runs and $m$ factors, where each factor takes $s$ levels from the set $\mathbb{Z}_{s}$.
Such a design can be represented by an $n\times m$ matrix $D=(x_{ik})$, where $x_{ik}\in\mathbb{Z}_{s}$ is the entry in the $i$-th row and $k$-th column, $i=1,\ldots,n$, $k=1,\ldots,m$. It is called a U-type design if the columns are balanced, that is, each level in $\mathbb{Z}_{s}$ appears equally often.

An orthogonal array (OA) of strength $t$,
denoted by $\mathrm{OA}(n,m,s_1\times\cdots\times s_m,t)$, is an $n\times m$ matrix whose $j$th column
takes values in $\mathbb{Z}_{s_j}$ and such that, in any of $t$-column subarray, all possible level
combinations appear equally often. When $s_1=\cdots=s_m=s$, we simply write $\mathrm{OA}(n,m,s,t)$.
A U-type  $(n, m, s)$ design is essentially an $\mathrm{OA}(n,m,s,1)$  \citep{Hedayat99}.
A Latin hypercube design (LHD) with $n$ runs and $m$ factors is a U-type  $(n, m, n)$ design, i.e., an $\mathrm{OA}(n,m,n,1)$ \citep{McKay79}.

A general strong orthogonal array (GSOA) with $n$ runs, $m$ factors, $s^p$ levels and strength $t$,
denoted by $\mathrm{GSOA}(n,m,s^p,t)$, is an $(n,m,s^p)$ design with entries from $\mathbb{Z}_{s^p}$ such that
any subarray of $g$ columns ($1\le g\le t$) can be collapsed into an
$\mathrm{OA}(n,g,s^{u_1}\times\cdots\times s^{u_g},g)$ for any positive integers
$\{u_1,\ldots,u_g\}$ satisfying $u_1+\cdots+u_g=t$ and $u_i\le p$ \citep{HeTang13, TianXu22}.
The collapsing map from $s^p$ levels to $s^{u}$ levels is
$\lfloor x/s^{p-u}\rfloor$ for $x\in\mathbb{Z}_{s^p}$, where $\lfloor\cdot\rfloor$ is the floor function.

\subsection{Stratified regions and the stratified \texorpdfstring{$L_2$}{L2}-discrepancy}

The stratified $L_2$-discrepancy is a new measure proposed by \cite{TianXu25} for evaluating the space-filling properties of point sets in stratified design regions.

Fix integers $s\ge 2$ and $p\ge 1$. For $u\in \mathbb{Z}_{p+1}= \{0,1,\ldots,p\}$ and $x\in[0,1)$, define the
\emph{stratified interval}
$$
R_u(x)=\Big[\frac{\lfloor s^u x\rfloor}{s^u},\ \frac{\lfloor s^u x\rfloor}{s^u}+\frac{1}{s^u}\Big),
$$
i.e., the unique interval among the $s^u$ equal subintervals of $[0,1)$ that contains $x$.
Obviously, $R_0(x) = [0, 1)$.
For $\bu=(u_1,\ldots,u_m)\in\mathbb{Z}_{p+1}^m$ and $\bx=(x_1,\ldots,x_m)\in[0,1)^m$, define the
$m$-dimensional stratified region as the Cartesian product
$
R_{\bu}(\bx)=R_{u_1}(x_1)\times\cdots\times R_{u_m}(x_m).
$
Its volume is $\text{Vol} (R_{\bu}(\bx)) =s^{-\sum_{j=1}^m u_j}$, which does not depend on $\bx$.

Suppose $P$ is a finite set of $n$ points in $[0,1)^m$.
The \emph{local projection discrepancy} of $P$ with respect to the stratification vector $\bu$ is then defined as
$$
\text{disc}^R_\bu(\bx; P) = \text{Vol}(R_\bu(\bx)) - \frac{|P\cap R_\bu(\bx)|}{n},
$$
where $\text{Vol}(R_\bu(\bx))$ denotes the volume of the region $R_\bu(\bx)$, and $|P \cap R_\bu(\bx)|$ represents the number of design points in $P$ that fall in $R_\bu(\bx)$.
Define the squared $L_2$-norm of $\text{disc}^R_\bu(\bx; P)$ as
$$
\left\|\text{disc}^R_\bu( P) \right\|^2=\int_{[0,1)^m}\left|\text{disc}^R_\bu(\bx; P)\right|^2\rd\bx.
$$

Let $\omega(\cdot)$ be a nonnegative weight function on $\mathbb{Z}_{p+1}$, and define the product-type weight
$
\omega(\bu)=\prod_{j=1}^m \omega(u_j)$, $\bu\in \mathbb{Z}_{p+1}^m,
$
with $\omega(\bu)>0$ for all $\bu$ of interest.
Here, weights associated with $\bu$ are incorporated in the definition of the stratified $L_2$-discrepancy to allow different stratification levels to be assigned different degrees of importance.
Then, the \emph{stratified $L_2$-discrepancy} is defined as
\[
\SD(P)=\left(\sum_{\bu \in \mathbb{Z}_{p+1}^m} \omega(\bu)\left\|\text{disc}^R_\bu( P) \right\|^2\right)^{\frac{1}{2}}.
\]
In this paper we will work with the \emph{squared} stratified $L_2$-discrepancy
\bea\label{SL2D}
\SD^2(P) = \sum_{\bu \in \mathbb{Z}_{p+1}^m} \omega(\bu)\, \left\|\text{disc}^R_\bu( P) \right\|^2.
\eea

As justified by \cite{TianXu25},  a suitable choice of the weight function follows an exponential weighting scheme
$\omega(\bu)=\prod_{j=1}^m \omega(u_j)=y^{\sum_{j=1}^m u_j}$,
with $y\in(0,1]$ being a constant parameter and $\omega(u_j) = y^{u_j}$.
Under this specification, the weights are proportional to the volumes of the corresponding stratified subregions, and therefore naturally reflect the space-filling hierarchy principle introduced by \cite{TianXu22}.
The parameter $y$ controls the relative emphasis placed on different stratification levels: when $y=1$ we have $\omega(\bu) = 1$ for all $\bu \in \mathbb{Z}_{p+1}^m$, which corresponds to constant weights, while smaller values of $y$ progressively downweight finer stratifications and favor uniformity at coarser resolutions.

Many $L_2$-type discrepancies admit a unified expression through a reproducing kernel Hilbert space (RKHS).
Let $\mathcal{X}=[0,1)^m$. Given a symmetric and nonnegative definite function $K: \mathcal{X}\times \mathcal{X} \to \mathbb{R}$, there exists a RKHS $\mathcal{H}$ of real-valued functions with the reproducing kernel $K$, i.e.,
$\mathcal{H}=\left\{F: \int_{\mathcal{X} \times \mathcal{X}}K(\bx,\by)\rd F(\bx)\rd F(\by)<\infty\right\}$,
and the inner product between two functions $F$ and $G\in\mathcal{H}$ is given by
$\langle F,G\rangle_{\mathcal{H}}=\int_{\mathcal{X} \times \mathcal{X}}K(\bx,\by)\rd F(\bx)\rd G(\by)$ \citep{Fang18}.
The norm on the space  $\mathcal{H}$ is then defined as $\|\cdot\|^2_{\mathcal{H}}=\langle \cdot, \cdot \rangle_{\mathcal{H}}$.

Given a set $P = \{\bz_1, \dots, \bz_n\}$ of $n$ points in the space $\mathcal{X} = [0,1)^m$ with $\bz_i=(z_{i1},\ldots,z_{im})$, let $F$ denote the cumulative distribution function of the uniform distribution on $[0,1)^m$ and $F_P$ be the empirical distribution induced by $P$.
The $L_2$-type discrepancy of $P$ associated with the reproducing kernel $K$ can be defined as the norm of the difference between $F$ and $F_P$ in the RKHS $\mathcal{H}$, i.e.,
\begin{align}\label{DPK}
{\rm Disc}(P;K)
&=\|F-F_P\|^2_{\mathcal{H}} \nonumber\\
&=\int_{\mathcal{X}\times\mathcal{X}} K(\bx,\by)\rd F(\bx)\rd F(\by)
-\frac{2}{n}\sum_{i=1}^n\int_{\mathcal{X}} K(\bz_i,\by)\rd F(\by) \nonumber\\
& ~~~ +\frac{1}{n^2}\sum_{i,j=1}^n K(\bz_i,\bz_j).
\end{align}
When the kernel is of product form
$
K(\bx,\by)=\prod_{k=1}^m f(x_k,y_k)$, $\bx,\by\in[0,1)^m$,
Equation \eqref{DPK} reduces to the classical closed form (see, e.g., \citealp{Hickernell98a,ZhouXu14})
\begin{eqnarray}\label{Disc}
{\rm Disc}(P;K)
&=&\left(\int_{[0,1)^2} f(x,y)\rd x\rd y\right)^{m}
-\frac{2}{n}\sum_{i=1}^{n}\prod_{k=1}^{m}\left(\int_{0}^{1} f(z_{ik},y)\rd y\right) \nonumber \\
&& +\frac{1}{n^2}\sum_{i,j=1}^{n}\prod_{k=1}^{m} f(z_{ik},z_{jk}).
\end{eqnarray}

For the stratified $L_2$-discrepancy \eqref{SL2D}, the associated kernel
\begin{equation}\label{KR}
K^{R}(\bx,\by)=\prod_{k=1}^{m} f^{R}(x_k,y_k),
\qquad
f^{R}(x,y)=\sum_{i=0}^{p}\frac{\omega(i)}{s^{i}}\delta_i(x,y),
\end{equation}
where, for $i=0,\ldots,p$,
\begin{equation}\label{Dfdi}
\delta_i(x,y)=\delta\!\left(\lfloor s^{i}x\rfloor,\ \lfloor s^{i}y\rfloor\right),
\end{equation}
and $\delta(\cdot,\cdot)$ is the Kronecker delta, i.e., $\delta(x,y)=1$ if $x=y$, and $0$ otherwise.
We see that $\delta_i(x,y)=1$ if and only if
$x$ and $y$ fall into the same subinterval when $[0,1)$ is partitioned into $s^i$ equal parts.
Substituting \eqref{KR} and \eqref{Dfdi} into \eqref{Disc} yields
\begin{eqnarray}\label{SD2}
\SD^2(P)
&=&\left(\int_{[0,1)^2} f^R(x,y)\rd x\rd y\right)^{m}
-\frac{2}{n}\sum_{i=1}^{n}\prod_{k=1}^{m}\left(\int_{0}^{1} f^R(z_{ik},y)\rd y\right) \nonumber \\
& & +\frac{1}{n^2}\sum_{i,j=1}^{n}\prod_{k=1}^{m} f^R(z_{ik},z_{jk})\nonumber\\
&=&-\left(\sum_{i=0}^{p}\frac{\omega(i)}{s^{2i}}\right)^{m}
+\frac{1}{n^2}\sum_{a,b=1}^{n}\prod_{k=1}^{m}
\left(\sum_{i=0}^{p}\frac{\omega(i)}{s^{i}}\delta_i(z_{ak},z_{bk})\right),
\end{eqnarray}
which coincides with Theorem~1 of \citet{TianXu25}.

\section{Uniform projection criterion under the stratified \texorpdfstring{$L_2$}{L2}-discrepancy}\label{sec3}

The centered $L_2$-discrepancy (CD) is one of the most widely used $L_2$-type measures of design uniformity; see, e.g., \citet{Fang06, Fang18}.
A well-known limitation of CD is its dimensionality effect, i.e., designs that appear highly uniform in the full space may exhibit poor uniformity in low-dimensional projections.
This issue is particularly relevant in computer experiments, where the number of inputs can be large and the response is often driven by a small subset of active variables according to the effect sparsity principle.
Motivated by these considerations, \citet{SunWangXu19} introduced the \emph{uniform projection criterion}, which averages the CD over all two-dimensional projections and seeks designs that minimize this average.
\citet{Liu23}  generalized the uniform projection criterion to an arbitrary $L_2$-type discrepancy.
Let $D=(x_{ik})$ be an $(n,m,s^p)$ design with entries in $\mathbb{Z}_{s^p}$.
Following the standard embedding for discrete designs, map $D$ to $P=\{ \bz_i\}_{i=1}^n =(z_{ik})$ in $[0,1)^m$ by
$
z_{ik}=({2x_{ik}+1})/{2s^p}$, $i=1,\ldots,n, k=1,\ldots,m.$
For a reproducing kernel $K$ on $[0,1)^m$, the $L_2$-type discrepancy of $D$ is defined by
${\rm Disc}(D,K)={\rm Disc}(P,K)$, where ${\rm Disc}(P,K)$ is given in~\eqref{Disc}.

The stratified $L_2$-discrepancy \eqref{SL2D} proposed by \citet{TianXu25} extends the generalized $L_2$ framework to settings in which the design space is stratified into subregions, thereby encouraging uniformity across multiple resolutions.
In many applications, however, projection behavior is more important than the full-dimensional space-filling, especially in screening stages when the design has large factor-to-run ratios \citep{Moon12,Joseph15,SunWangXu19}.
Therefore, we study the uniform projection criterion induced by the stratified $L_2$-discrepancy.
Throughout this paper, $\Phi_{\rm SD}(D)$ denotes the average squared stratified $L_2$-discrepancy over all two-dimensional projections,
\begin{equation}\label{eqphiDisc}
\Phi_{\rm SD}(D)=\frac{2}{m(m-1)}\sum_{|u|=2} \SD^2(D_u),
\end{equation}
where $u\subset\{1,\ldots,m\}$ indexes a two-factor projection, $|u|=2$, and $D_u$ denotes the subdesign of $D$ formed by the columns in $u$.
A design $D$ is said to be a \emph{uniform projection design} under the stratified $L_2$-discrepancy if it has the minimum  $\Phi_{\rm SD}(D)$ value among all $(n,m,s^p)$ designs.

\begin{defi}\label{def:habdist}
For an $(n,m,s^p)$ design $D=(x_{ik})$ and $a,b\in\{1,\ldots,n\}$, define the \emph{weighted hierarchical distance}  between the $a$-th and $b$-th runs of $D$ as
\begin{equation}\label{eq:dab_new}
\dist_{ab}
:=\sum_{k=1}^m\sum_{i=0}^p \frac{\omega(i)}{s^i}\Bigl\{1-\delta_i\bigl(z_{ak},z_{bk}\bigr)\Bigr\},
\end{equation}
where $z_{ik}=(2x_{ik}+1)/(2s^p)$.
\end{defi}

We assume that $\omega(i)\ge 0$ for all $i$ and $\omega(i)>0$ for at least one $i\ge 0$.
The quantity $\dist_{ab}$ is a metric on the point set $P=\{ \bz_a\}_{a=1}^n$ and $\dist_{ab}=0$ implies $\bz_a=\bz_b$.
It can be viewed as a weighted multi-resolution Hamming distance, i.e., a mismatch at resolution $i$ contributes $\omega(i)/s^i$ in each coordinate.

The following result gives an analytical expression of $\Phi_{\rm SD}$ defined in \eqref{eqphiDisc}.

\begin{theorem}\label{thmSD2}
Let $D=(x_{ik})$ be a U-type $(n,m,s^p)$ design and let $z_{ik}=(2x_{ik}+1)/(2s^p)$.
Then we have
\begin{align*}
\Phi_{\rm SD}(D)
&=\frac{1}{n^2m(m-1)}\sum_{a=1}^n\sum_{b=1}^n \dist_{ab}^2
+ C_{\rm SD}(m,s,p),
\end{align*}
where the constant $C_{\rm SD}(m,s,p)$ depends only on $(m,s,p)$ and the weight function $\omega(\cdot)$, and is given by
\begin{align}\label{eq:CSD_simplified}
C_{\rm SD}(m,s,p)
&= -\frac{m}{m-1}A_0^2
+\frac{2m}{m-1}A_0A_1
-A_1^2
-\frac{1}{m-1}B
-\frac{2}{m-1}C
\end{align}
with
\[
A_0=\sum_{i=0}^p\frac{\omega(i)}{s^i},\qquad
A_1=\sum_{i=0}^p\frac{\omega(i)}{s^{2i}},\qquad
B=\sum_{i=0}^p\frac{\omega(i)^2}{s^{3i}},\qquad
C=\sum_{0\le i<j\le p}\frac{\omega(i)\omega(j)}{s^{i+2j}}.
\]
\end{theorem}

By the definition \eqref{eqphiDisc},  to obtain a design's $\Phi_{\rm SD}$ value, we first calculate the discrepancies of all the column pairwise projections of $D$ and then take an average. Therefore, \eqref{eqphiDisc} can be viewed as a measure of a design's column pairwise relationships.
Theorem~\ref{thmSD2} establishes a new formula for \eqref{eqphiDisc} for U-type designs from the viewpoint of row-pairwise weighted hierarchical distance \eqref{eq:dab_new}.

\begin{example}\label{ex1}
  As a specific case of Theorem~\ref{thmSD2}, we consider the use of exponential weights $\omega(i)=y^i$ with $y\in(0,1]$. Let $$S_p(q) = \sum_{i=0}^p q^i = \frac{1-q^{p+1}}{1-q} \text{ for } q \in (0,1).$$ Under this weighting scheme, after some tedious calculations, the projection discrepancy $\Phi_{\rm SD}(D)$ can be expressed as
\begin{equation*}
\Phi_{\rm SD}(D) = \frac{1}{n^2m(m-1)}\sum_{a=1}^n\sum_{b=1}^n \dist_{ab}^2 + C_{\rm SD}^{(y)}(m,s,p),
\end{equation*}
where the weighted hierarchical distance $\dist_{ab}$ for points $a,b \in \{1,\dots,n\}$ is defined as
\begin{equation*}
\dist_{ab} = \sum_{k=1}^m \sum_{i=0}^p \frac{y^i}{s^i} \bigl\{1-\delta_i(z_{ak},z_{bk})\bigr\},
\end{equation*}
and the constant $C_{\rm SD}^{(y)}(m,s,p)$ depends solely on $(m,s,p,y)$ and is given by
\begin{equation*}
{C_{\rm SD}^{(y)}(m,s,p) = -\frac{m}{m-1}A_0^2 + \frac{2m}{m-1}A_0A_1 - A_1^2 - \frac{1}{m-1}B - \frac{2}{m-1}C,}
\end{equation*}
with
\begin{align*}
{A_0 = S_p\left(\frac{y}{s}\right), \quad A_1 = S_p\left(\frac{y}{s^2}\right), \quad B = S_p\left(\frac{y^2}{s^3}\right), \quad C = \frac{S_p(y/s^2) - S_p(y^2/s^3)}{1 - y/s}.}
\end{align*}

\end{example}

Now consider lower and upper bounds for the uniform projection criterion $\Phi_{\rm SD}(D)$.
It suffices to bound
\[
G_D:=\sum_{a=1}^n\sum_{b=1}^n \dist_{ab}^2,
\]
where $\dist_{ab}$ is the weighted hierarchical distance defined in \eqref{eq:dab_new},
since
\begin{equation}\label{eq:PhiSD_GD}
\Phi_{\rm SD}(D)=\frac{G_D}{n^2m(m-1)}+C_{\rm SD}(m,s,p)
\end{equation}
by Theorem~\ref{thmSD2}.
The next result provides bounds for $G_D$, and hence for $\Phi_{\rm SD}(D)$.

\begin{theorem}\label{thm:bounds_PhiSD}
Let $D$ be a U-type $(n,m,s^p)$ design. Let $A_0$, $A_1$, $B$ and $C$ be the constants defined in Theorem~\ref{thmSD2}, and define
$
A_0^{(\ell)}:=\sum_{i=0}^{\ell}\frac{\omega(i)}{s^i}
$  for $\ell=0,\ldots,p-1$.
Then we have
$
G_D^{\rm LB}\le G_D\le G_D^{\rm UB},
$
where
$$
G_D:=\sum_{a=1}^n\sum_{b=1}^n \dist_{ab}^2,
\qquad
G_D^{\rm LB}
=\frac{n^3m^2}{n-1}\bigl(A_0-A_1\bigr)^2,
\qquad
G_D^{\rm UB}
=n^2m^2\sum_{\ell=0}^{p-1}\frac{s-1}{s^{\ell+1}}\Bigl(A_0-A_0^{(\ell)}\Bigr)^2.
$$
Moreover, the lower bound $G_D^{\rm LB}$ is attained if and only if $\dist_{ab}=\bar d$ for all $1\le a<b\le n$, where
$$
\bar d=\frac{nm}{n-1}\bigl(A_0-A_1\bigr).
$$
The upper bound $G_D^{\rm UB}$ is attained by any U-type $(n,m,s^p)$ design whose $m$ columns are identical copies of a balanced column on $\mathbb{Z}_{s^p}$, i.e., a column in which each of the $s^p$ levels appears $n/s^p$ times.
Consequently, by \eqref{eq:PhiSD_GD} we have
$$
\Phi_{\rm SD}^{\rm LB}\le \Phi_{\rm SD}(D)\le \Phi_{\rm SD}^{\rm UB},
$$
where
\begin{align}\label{eq:Phi_LB_collected_inthm}
\Phi_{\rm SD}^{\rm LB}
&=\frac{m}{(n-1)(m-1)}\,A_0^2
-\frac{2m}{(n-1)(m-1)}\,A_0A_1
+\frac{n+m-1}{(n-1)(m-1)}\,A_1^2
-\frac{1}{m-1}\,B-\frac{2}{m-1}\,C,
\end{align}
and
\begin{align}\label{eq:Phi_UB_collected_inthm_noTp}
\Phi_{\rm SD}^{\rm UB}
&=\frac{m}{(m-1)s^p}A_0^2
-A_1^2
+\frac{m}{m-1}\sum_{\ell=0}^{p-1}\frac{s-1}{s^{\ell+1}}\Bigl(A_0^{(\ell)}\Bigr)^2 -\frac{1}{m-1}\,B-\frac{2}{m-1}\,C.
\end{align}
\end{theorem}

\begin{example}\label{ex2}
Continuing Example~\ref{ex1}, we specialize Theorem~\ref{thm:bounds_PhiSD} to the same exponential weights
$\omega(i)=y^i$ with $y\in(0,1]$.  Recall the shorthand $S_p(q)$  introduced in Example~\ref{ex1}.
Then
\[
A_0^{(\ell)}=S_{\ell}\!\left(\frac{y}{s}\right)\ \ (\ell=0,\ldots,p-1),
\]
and $A_0$, $A_1$, $B$ and $C$ are given in Example~\ref{ex1}. Substituting these quantities into
Theorem~\ref{thm:bounds_PhiSD} yields
\[
G_D^{\rm LB}
=\frac{n^3m^2}{n-1}\left\{S_p\!\left(\frac{y}{s}\right)-S_p\!\left(\frac{y}{s^2}\right)\right\}^2,
\]
and
\[
G_D^{\rm UB}
=n^2m^2\sum_{\ell=0}^{p-1}\frac{s-1}{s^{\ell+1}}
\left\{S_p\!\left(\frac{y}{s}\right)-S_{\ell}\!\left(\frac{y}{s}\right)\right\}^2.
\]
Moreover, the equality condition for the lower bound becomes $\dist_{ab}=\bar d$ for all $1\le a<b\le n$, where
\[
\bar d=\frac{nm}{n-1}\left\{S_p\!\left(\frac{y}{s}\right)-S_p\!\left(\frac{y}{s^2}\right)\right\}.
\]
Finally, the bounds $\Phi_{\rm SD}^{\rm LB}$ and $\Phi_{\rm SD}^{\rm UB}$ in
\eqref{eq:Phi_LB_collected_inthm}--\eqref{eq:Phi_UB_collected_inthm_noTp}
follow by plugging the same specialization of $A_0,A_1,A_0^{(\ell)},B$ and $C$ into those formulas.
\end{example}

\begin{remark}
For $m\geq 3$, define the average squared stratified $L_2$-discrepancy over all three-column projections by
\[
\Phi_{\rm SD,3}(D) = \frac{6}{m(m-1)(m-2)}\sum_{|u|=3}\mathrm{SD}^2(D_u),
\]
which is analogous to the average centered $L_2$-discrepancy over all three-dimensional projections in \citet{Zhou23phd}. When $m=3$, this criterion reduces to the ordinary squared stratified $L_2$-discrepancy.

For a general U-type design, $\Phi_{\rm SD,3}(D)$ cannot be expressed only through the pairwise distances $d_{ab}$ and constants. It also contains the interaction term $\sum_{a,b} d_{ab} q_{ab}$, where $q_{ab} = \sum_{k=1}^m \left(\sum_{i=0}^p {\omega(i)}{s^{-i}}\left[1-\delta_i\bigl(z_{ak},z_{bk}\bigr)\right]\right)^2$,
which makes the general three-dimensional problem difficult.

A simpler conclusion can be obtained within the class of orthogonal arrays of strength two. If $D$ is an $\mathrm{OA}(n,m,s^p,2)$, then the term $\sum_{a,b}d_{ab}q_{ab}$ is constant. Hence, up to constants, the design-dependent part of $\Phi_{\rm SD,3}(D)$ becomes
\[
\frac{1}{n^2m(m-1)(m-2)}
\sum_{a,b=1}^n
\left\{3A_0(m-2)d_{ab}^2-d_{ab}^3\right\},
\]
where $A_0=\sum_{i=0}^p \omega(i)/s^i$.
For $m \geq 2A_0/\omega(0)$, the function $f(x)=3A_0(m-2)x^2-x^3$ is strictly convex over the admissible range of $d_{ab}$. Consequently, if there exists an OA$(n,m,s^p,2)$ such that all $d_{ab}$ are equal, it attains the minimum value of $\Phi_{\rm SD,3}(D)$ within the class of strength-two OAs.
Specifically,  when $\omega(i)\in(0, 1]$ for all $i$, it follows that $A_0\leq \sum_{i=0}^p s^{-i}\leq 2$. Hence the condition $m\omega(0)\geq 2A_0$ holds for all $m\geq 4$. The existence of such equal-distance OAs and the general U-type case without the strength-two OA restriction deserve further study.
\end{remark}

\section{Optimal space-filling designs under the \texorpdfstring{$\Phi_{\rm SD}$}{Phi SD} criterion}\label{sec3opt}

Theorem~\ref{thmSD2} shows that, for any U-type $(n,m,s^p)$ design $D$, minimizing $\Phi_{\rm SD}(D)$ over U-type designs is equivalent to minimizing $G_D=\sum_{a=1}^n\sum_{b=1}^n \dist_{ab}^2$.
Theorem~\ref{thm:bounds_PhiSD} provides explicit lower and upper bounds for $G_D$ (and hence for $\Phi_{\rm SD}(D)$). In particular, the lower bound is attained if and only if all pairwise distances are equal, namely
\[
\dist_{ab}=\bar d \quad\text{for all }1\le a<b\le n,
\qquad
\bar d=\frac{nm}{n-1}\bigl(A_0-A_1\bigr) = \sum_{i=0}^{p} \omega(i) \frac{(s^i-1)nm}{s^{2i}(n-1)}.
\]
A convenient sufficient condition that guarantees $\dist_{ab}\equiv \bar d$ \emph{regardless of the particular choice of weights} is the following balance condition at every resolution:
\bea\label{dtave}
\sum_{k=1}^m\delta_i(z_{ak},z_{bk})
=\frac{nm}{s^i(n-1)}-\frac{m}{n-1}
\quad\text{for all }1\le a<b\le n,\ \ 0\le i\le p.
\eea
Indeed, \eqref{dtave} fixes the number of coordinates that coincide when $[0,1)$ is partitioned into $s^i$ equal subintervals for each pair of runs $(a,b)$. Equivalently, for each $i=0,\ldots,p-1$, exactly
$
{nm(s-1)}/[{s^{i+1}(n-1)}]
$
of the $m$ coordinate-pairs $\{z_{ak},z_{bk}\}$ (with $k=1,\ldots,m$) lie in the same $s^i$-subinterval but split when the partition is refined to $s^{i+1}$ subintervals, and exactly
$
{nm}/[{s^p(n-1)}]-{m}/({n-1})
$
of them remain in the same subinterval under the finest partition into $s^p$ equal parts. Moreover, if the weights are allowed to take arbitrary values, then attaining the lower bound (uniform pairwise distances) forces these counts to be the same for every pair, and hence \eqref{dtave} is also necessary.

In the following, we present several constructions of optimal designs that satisfy \eqref{dtave}; these designs therefore attain the lower bounds $G_D^{\rm LB}$ and $\Phi_{\rm SD}^{\rm LB}$ in Theorem~\ref{thm:bounds_PhiSD}.

\subsection{Construction of optimal designs using the multiplication table of a Galois field}

For any $y\in \mathbb{Z}_{s^p}$, we can write
\begin{equation}\label{numdig}
y=\sum_{j=1}^p f_j(y)\,s^{p-j},
\end{equation}
where $f_j(y)\in\mathbb{Z}_s$ is the $j$th digit of $y$ in its base-$s$ expansion.
The (Niederreiter-Rosenbloom-Tsfasman) NRT-distance \citep{Bierbrauer2002} between  $x$ and $y\in\mathbb{Z}_{s^p}$ is defined by
\begin{equation}\label{eq:NRT-dist}
\rho(x,y)=
\begin{cases}
p+1-\min\{\,j: f_j(x)\neq f_j(y),\ j=1,\ldots,p\,\}, & x\neq y,\\
0, & x=y.
\end{cases}
\end{equation}
 When $p=1$, the NRT-distance is the Hamming distance, that is,
$\rho(x,y)=1$ if $x\neq y$ and $0$ otherwise.

To describe algebraic constructions that attain the lower bounds in Theorem~\ref{thm:bounds_PhiSD}, we briefly recall some basic concepts of the Galois field.  Let $\mathrm{GF}(s^p)$ be the Galois field of order $s^p$, where $s$ is prime.  Its elements can be represented as polynomials over $\mathbb{Z}_s$ of degree less than $p$.
For
$
y=\sum_{j=1}^p c_j \xi^{p-j}\in \mathrm{GF}(s^p),
$
where $c_j\in\mathbb{Z}_s$ and $\xi$ is an indeterminate variable, define the coordinate map
$f_j:\mathrm{GF}(s^p)\to\mathbb{Z}_s$ by $f_j(y)=c_j$ for $j=1,\ldots,p$.
This mapping extracts the coefficient of $\xi^{p-j}$ in $y$.
Under the standard identification between $\mathrm{GF}(s^p)$ and $\mathbb{Z}_{s^p}$ via the digit vector $(f_1(y),\ldots,f_p(y))$, the number representation in \eqref{numdig} is consistent with this coordinate system. We also write $\sum_{j=1}^p f_j(y) s^{p-j}$  as the number representation of $y \in \mathrm{GF}(s^p)$.

Addition and subtraction in $\mathrm{GF}(s^p)$ are the usual polynomial operations over $\mathbb{Z}_s$, while multiplication is taken modulo an irreducible polynomial of degree $p$.  The addition table of $\mathrm{GF}(s^p)$ is uniquely determined, but different choices of the irreducible polynomial may lead to different multiplication tables of $\mathrm{GF}(s^p)$, and hence to different designs.

Let $M$ be the multiplication table of $\mathrm{GF}(s^p)$.  By deleting its first column (corresponding to multiplication by $0$), \citet[Theorem~6]{TianXu24} constructed a design with $s^p$ rows and $s^p-1$ columns (more precisely, a $\mathrm{GSOA}(s^p,s^p-1,s^p,t)$, see Proposition~\ref*{CorTX}~(i) below for more details) in which, for any two distinct rows, the numbers of coordinatewise pairs at each NRT-distance are the same (i.e., the NRT-distance profile is constant over all row pairs).  This property implies the pair-count condition \eqref{dtave}, which is the equality condition for the lower bound in Theorem~\ref{thm:bounds_PhiSD}.  Consequently, the resulting design attains $G_D^{\rm LB}$ and therefore also achieves $\Phi_{\rm SD}^{\rm LB}$ via \eqref{eq:PhiSD_GD}.  Furthermore, when $s$ is an odd prime, an optimal $\mathrm{GSOA}(s^p,(s^p-1)/2,s^p,t)$  can be obtained by selecting half of the columns from a $\mathrm{GSOA}(s^p,s^p-1,s^p,t)$ derived from the full multiplication table. Note that these designs also attain the lower bound for the stratified $L_2$-discrepancy SD$^2(D)$ in \citet[Theorem~5]{TianXu25}.

More generally, for any positive integer $q < p$, one can collapse the~$s^p$ levels to~$s^q$ levels by truncating the last $p-q$ base-$s$ digits: map
$
y=\sum_{j=1}^p f_j(y)\,\xi^{p-j}
$
to
$
y'=\sum_{j=1}^q f_j(y)\,\xi^{q-j}.
$
Deleting the first column of the resulting (collapsed) multiplication table again yields a design (a GSOA, see Proposition~\ref*{CorTX}~(i)) satisfying \eqref{dtave}, and hence a $\Phi_{\rm SD}$-optimal design under Theorem~\ref{thm:bounds_PhiSD}.

The following existence results summarize the constructions in \citet{TianXu24,TianXu25}.

\begin{prop}\label{CorTX}
\begin{itemize}
  \item[(i)]
  For any prime number $s$ and any positive integers  $q\le p$, there exists a $\mathrm{GSOA}(s^p,\,s^p-1,\,s^q,\,t)$ attaining the lower bounds $G_D^{\rm LB}$ and $\Phi_{\rm SD}^{\rm LB}$ in Theorem~\ref{thm:bounds_PhiSD}, with $t=2$ if $s=2$ and $t=1$ if $s>2$.
  \item[(ii)]
  For any odd prime $s$, there exists a $\mathrm{GSOA}\bigl(s^p,\,(s^p-1)/2,\,s^p,\,t\bigr)$ attaining the lower bounds $G_D^{\rm LB}$ and $\Phi_{\rm SD}^{\rm LB}$ in Theorem~\ref{thm:bounds_PhiSD}, with $t=2$ if $s=3$ and $t=1$ if $s>3$.
\end{itemize}
\end{prop}

As illustrative examples of Proposition~\ref*{CorTX}, \citet{TianXu24} constructed a $\mathrm{GSOA}(16,15,2^4,2)$ from the multiplication table of $\mathrm{GF}(2^4)$, and a $\mathrm{GSOA}(16,15,2^3,2)$ from a collapsed multiplication table with $8$ levels.  They also generated families of nonisomorphic designs by adding field constants to the multiplication-table construction, and then selected efficient designs using a secondary criterion.

Consider case (ii) of Proposition~\ref*{CorTX}. \citet[Example~5]{TianXu25} constructed a $\mathrm{GSOA}(9,8,3^2,1)$ by deleting the first column of the multiplication table over $\mathrm{GF}(3^2)$.  Partitioning the nonzero elements of $\mathrm{GF}(3^2)$ into two disjoint sets $\{1,3,4,5\}$ and $\{2,6,8,7\}$, the subarray formed by the columns indexed by $\{1,3,4,5\}$ yields an optimal $\mathrm{GSOA}(9,4,3^2,2)$ with respect to $\Phi_{\rm SD}$.

For convenience, all the above mentioned designs can be found in the Supplementary Materials.

\subsection{More methods for constructing optimal designs}

\citet{Gao25} extended the construction framework of optimal space-filling designs developed by \citet{TianXu24,TianXu25}, thereby enlarging the parameter range and improving flexibility. As noted in the previous subsection, Tian and Xu constructed optimal designs from the multiplication table of a Galois field $\mathrm{GF}(s^p)$  (with $s$ prime).  Such a multiplication table is a special case of a generalized Hadamard matrix. Building on this connection, \citet{Gao25} developed generalized Hadamard matrix constructions for any prime power $s$ (see Proposition~\ref*{prop:GJW} (i) below for more details), which substantially broadens the class of admissible fields.   They also introduced a Rao--Hamming type construction based on saturated linear orthogonal arrays, leading to optimal designs such as $\mathrm{GSOA}(s^2,s+1,s^2,2)$.
As \citet{Gao25} is an unpublished manuscript, all design matrices cited from this work are provided in the Supplementary Materials.

The optimality results in \citet{Gao25} are stated under the minimum aberration-type space-filling criterion and are characterized by a uniformity condition on pairwise NRT-distance counts. For the designs produced by their constructions, this condition implies \eqref{dtave}; consequently, by Theorem~\ref{thm:bounds_PhiSD}, these designs attain the lower bounds $G_D^{\rm LB}$ and $\Phi_{\rm SD}^{\rm LB}$. The following Proposition~\ref*{prop:GJW} summarizes \citet{Gao25}'s constructions.

Let $s$ be a prime power and let $\mathrm{GF}(s)$ be the finite field of order $s$.
An $n\times n$ array with entries in $\mathrm{GF}(s)$ is called a \emph{generalized Hadamard matrix} \citep{Hedayat99}, denoted by \(\mathrm{GH}(n, s)\),  if for every pair of distinct columns,  each element of $\mathrm{GF}(s)$ appears exactly $n/s$ times in the componentwise difference of those columns.

\begin{prop}\label{prop:GJW}
\begin{itemize}
  \item[(i)] For any prime power $s$, and any positive integers $\lambda$ and $q\leq p$, suppose a generalized Hadamard matrix $\mathrm{GH}(\lambda s^p,$ $s^p)$ exists. Then, there exists a
  $\mathrm{GSOA}(\lambda s^p,\lambda s^p-1,s^q,t)$ that attains the lower bounds $G_D^{\rm LB}$ and $\Phi_{\rm SD}^{\rm LB}$ in Theorem~\ref{thm:bounds_PhiSD}, where $t=2$ if $s=2$ and $t=1$ if $s>2$.
  \item[(ii)]
  For any prime power $s$, there exists a balanced $\mathrm{GSOA}(s^2,s+1,s^2,2)$ that attains the lower bounds
  $G_D^{\rm LB}$ and $\Phi_{\rm SD}^{\rm LB}$ in Theorem~\ref{thm:bounds_PhiSD}.
\end{itemize}
\end{prop}

\citet{Gao25} constructed $\mathrm{GSOA}(32,31,4^2,1)$ and $\mathrm{GSOA}(32,31,4^1,1)$ (with $s=4$), {$\mathrm{GSOA}(32,31,16^1,1)$ (with $s=16$)}, and $\mathrm{GSOA}(9,4,3^2,2)$ as examples of optimal designs, all of which also illustrate  Proposition~\ref*{prop:GJW}.
They further proposed two operations, juxtaposition and collapsing, to generate new optimal designs from existing ones.

\smallskip
\noindent\textbf{Juxtaposition.}
Let $D_1$ and $D_2$ be balanced $\mathrm{GSOA}(n,m_1,s^p,t_1)$ and $\mathrm{GSOA}(n,m_2,s^p,t_2)$, respectively, and suppose that both satisfy \eqref{dtave}. Then the concatenated design
$D=(D_1,D_2)$ is a balanced $\mathrm{GSOA}(n,m_1+m_2,s^p,t)$ with $1\le t\le \min\{t_1,t_2\}$, and it also satisfies \eqref{dtave} (with $m=m_1+m_2$). Hence $D$ attains $G_D^{\rm LB}$ and $\Phi_{\rm SD}^{\rm LB}$ in Theorem~\ref{thm:bounds_PhiSD}.
For example, taking $D_1=\mathrm{GSOA}(9,8,3^2,1)$ and $D_2=\mathrm{GSOA}(9,4,3^2,2)$ from Example~3 of \citet{Gao25} yields an optimal design $\mathrm{GSOA}(9,12,3^2,1)$.

\smallskip
\noindent\textbf{Collapsing.}
Let $\varphi:\mathbb{Z}_{s^p}\to \mathbb{Z}_{s^q}$ ($q\le p$) be the digit-truncation map
\[
\varphi(z)=\sum_{i=1}^q f_i(z)\,s^{\,q-i},
\qquad
z=\sum_{i=1}^p f_i(z)\,s^{\,p-i}\in\mathbb{Z}_{s^p},
\]
where $f_i(\cdot)$ is defined in \eqref{numdig}. For a matrix $D=(d_{ij})$ with entries in $\mathbb{Z}_{s^p}$, define
$\varphi(D)=\bigl(\varphi(d_{ij})\bigr)$.
If $D$ is a balanced $\mathrm{GSOA}(n,m,s^p,t)$ satisfying \eqref{dtave}, then $\varphi(D)$ is a balanced $\mathrm{GSOA}(n,m,s^q,t)$ and satisfies \eqref{dtave} with $p$ replaced by $q$; therefore, $\varphi(D)$ also attains $G_D^{\rm LB}$ and $\Phi_{\rm SD}^{\rm LB}$ in Theorem~\ref{thm:bounds_PhiSD}.
As an illustration (Example~4 of \citealp{Gao25}), collapsing an optimal $\mathrm{GSOA}(9,12,3^2,1)$ produces an optimal $\mathrm{GSOA}(9,12,3^1,1)$ under the $\Phi_{\rm SD}$ criterion.

All the above mentioned designs can be found in the Supplementary Materials.

 \section{Numerical studies}\label{sec4}

This section verifies the properties of the proposed uniform projection criterion $\Phi_{\rm SD}$ through numerical examples. Following \citet{TianXu25}, we adopt the same kernel choices for reproducibility.

We first consider two GSOAs from  Example~5 of \citet{TianXu25}, $\mathrm{GSOA}(9,8,3^2,1)$ and $\mathrm{GSOA}(9,4,3^2,2)$, which are optimal under the SD.  These two designs can be found in the Supplementary Materials. According to Proposition~\ref{CorTX}, these designs are also optimal under $\Phi_{\rm SD}$, attaining the lower bounds $\Phi_{\rm SD}^{\rm LB}$ and $G_D^{\rm LB}$ established in Theorem~\ref{thm:bounds_PhiSD}.

We recompute the $\SD^2$ values reported in Example~5 of \citet{TianXu25}, which were computed under constant weights. Furthermore, we evaluate $\Phi_{\rm SD}$ and $G_D$, alongside their theoretical bounds for the two designs. Table~\ref{tab:sec4_ex5_sd2_updsd} reports the computed values. The numerical results confirm that the criteria values for these optimal designs exactly match their respective lower bounds, verifying the sharpness of the bounds derived in Section~\ref{sec3}.

To provide a simple performance range, we further compare each optimal design with (i) a worst-case U-type design (all columns identical) and (ii) 100 random U-type designs. The worst-case designs reach the upper bound $\Phi_{\rm SD}^{\rm UB}$. For the $(9,8,3^2,1)$ case, the random designs yield a mean $\Phi_{\rm SD}$ of $0.012899$ (sd: $0.000773$). For the $(9,4,3^2,2)$ case, the mean is $0.012888$ (sd: $0.001615$). These results demonstrate that the optimal GSOAs consistently reside at the theoretical minimum of the $\Phi_{\rm SD}$ range.

\begin{table}[htbp]
\centering
\caption{Criteria values and theoretical bounds for optimal GSOAs.}
\label{tab:sec4_ex5_sd2_updsd}
\setlength{\tabcolsep}{4pt}
\begin{tabular}{l c cc cc cc}
\toprule
Design & $\SD^2$ & $\Phi_{\rm SD}$ & $\Phi_{\rm SD}^{\rm LB}$ & $\Phi_{\rm SD}^{\rm UB}$ & $G_D$ & $G_D^{\rm LB}$ & $G_D^{\rm UB}$ \\
\midrule
$\mathrm{GSOA}(9,8,3^2,1)$ & 1.148028 & 0.010234 & 0.010234 & 0.031398 & 600.888889 & 600.888889 & 696.888889 \\
$\mathrm{GSOA}(9,4,3^2,2)$ & 0.075833 & 0.006706 & 0.006706 & 0.031398 & 150.222222 & 150.222222 & 174.222222 \\
\bottomrule
\end{tabular}
\end{table}

Next, we consider a standard benchmark involving four $19 \times 18$ LHDs which have been studied by \cite{SunWangXu19} and \cite{TianXu25}. The four LHDs are: a maximin distance design \cite{Johnson90}, a maximum projection (MaxPro) design \cite{Joseph15}, a uniform design under the CD, and a uniform projection design under the CD (UPD). This setup is frequently used to assess low-dimensional projection uniformity, where UPD is expected to perform well by construction.

Following Section~4 of \citet{TianXu25},  we evaluate two kernel settings indexed by $s \in \{2, 3\}$, with the stratification parameter $p = \lfloor \log_s(n) \rfloor$ with $n=19$, and use constant weights for the four designs. The resulting $\Phi_{\rm SD}$ criteria are denoted by $\Phi_{\rm SD}^{(2)}$ and $\Phi_{\rm SD}^{(3)}$, respectively.

Table~\ref{tab:sec4_sdupd_19x18} augments the full-dimensional criteria (CD, WD, MD) and stratified discrepancy (SD) values from \citet{TianXu25} with our $\Phi_{\rm SD}$ criteria. Numerical results show that while CD, WD, and MD each select different designs as the ``best'', both the SD and our proposed $\Phi_{\rm SD}$ (under $s=2$ and $s=3$) consistently rank the UPD as the most space-filling design. This consistency arises because $\Phi_{\rm SD}$ is specifically designed to capture low-dimensional projection properties. Notably, UPD achieves the smallest $\Phi_{\rm SD}$ values across both kernel settings, indicating that our criterion places a stronger emphasis on projection-averaged stratified uniformity than the other full-dimensional discrepancies compared {here}.

\begin{table}[htbp]
\centering
\caption{Comparison of four $19 \times 18$ LHDs using full-dimensional and projection-based criteria. Bold values indicate the best design under each criterion.}
\label{tab:sec4_sdupd_19x18}
\setlength{\tabcolsep}{3pt}
\begin{tabular}{lccccccc}
\toprule
Design & CD & WD & MD & SD(\({s}=2\)) & SD(\({s}=3\)) & \(\Phi_{\rm SD}^{(2)}\) & \(\Phi_{\rm SD}^{(3)}\) \\
\midrule
Maximin & 1.2889 & 7.0488 & 25.2549 & 87.7170 & 6.0710 & 0.0227647 & 0.0072553\\
MaxPro  & 1.3090 & \textbf{6.8823} & 24.8515 & 87.6938 & 6.0468 & 0.0221945 & 0.0068359\\
Uniform & \textbf{1.2643} & 6.9414 & \textbf{24.8049} & 87.6903 & 6.0496 & 0.0220757 & 0.0068860\\
UPD     & 1.2655 & 6.9352 & 24.8554 & \textbf{87.6342} & \textbf{6.0365} & \textbf{0.0204921} & \textbf{0.0066817}\\
\bottomrule
\end{tabular}
\end{table}

\begin{figure}[htbp]
\centering
\begin{minipage}[t]{0.48\textwidth}
\centering
\includegraphics[width=\textwidth]{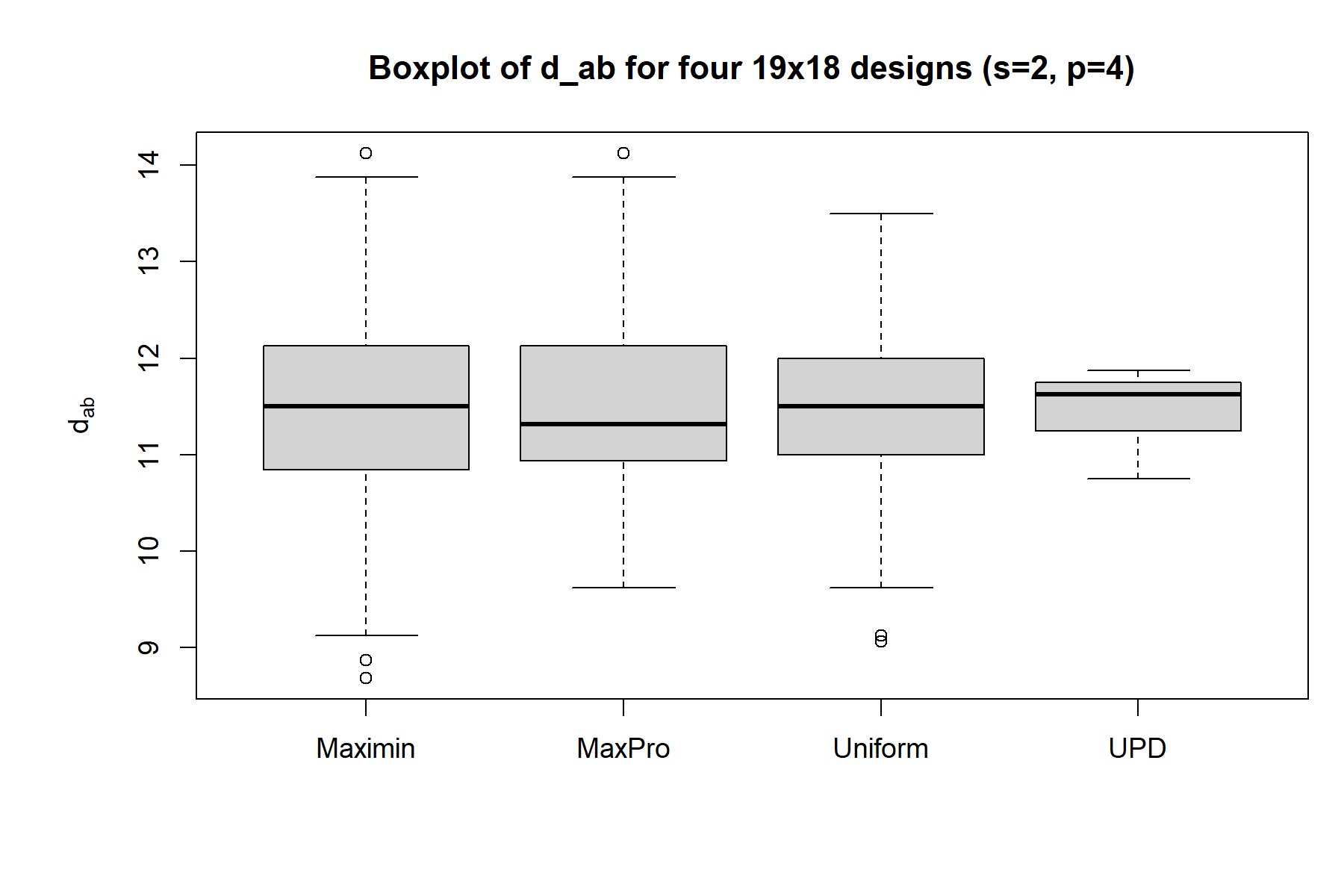}\\
\small (a) $s=2$, $p=4$
\end{minipage}\hfill
\begin{minipage}[t]{0.48\textwidth}
\centering
\includegraphics[width=\textwidth]{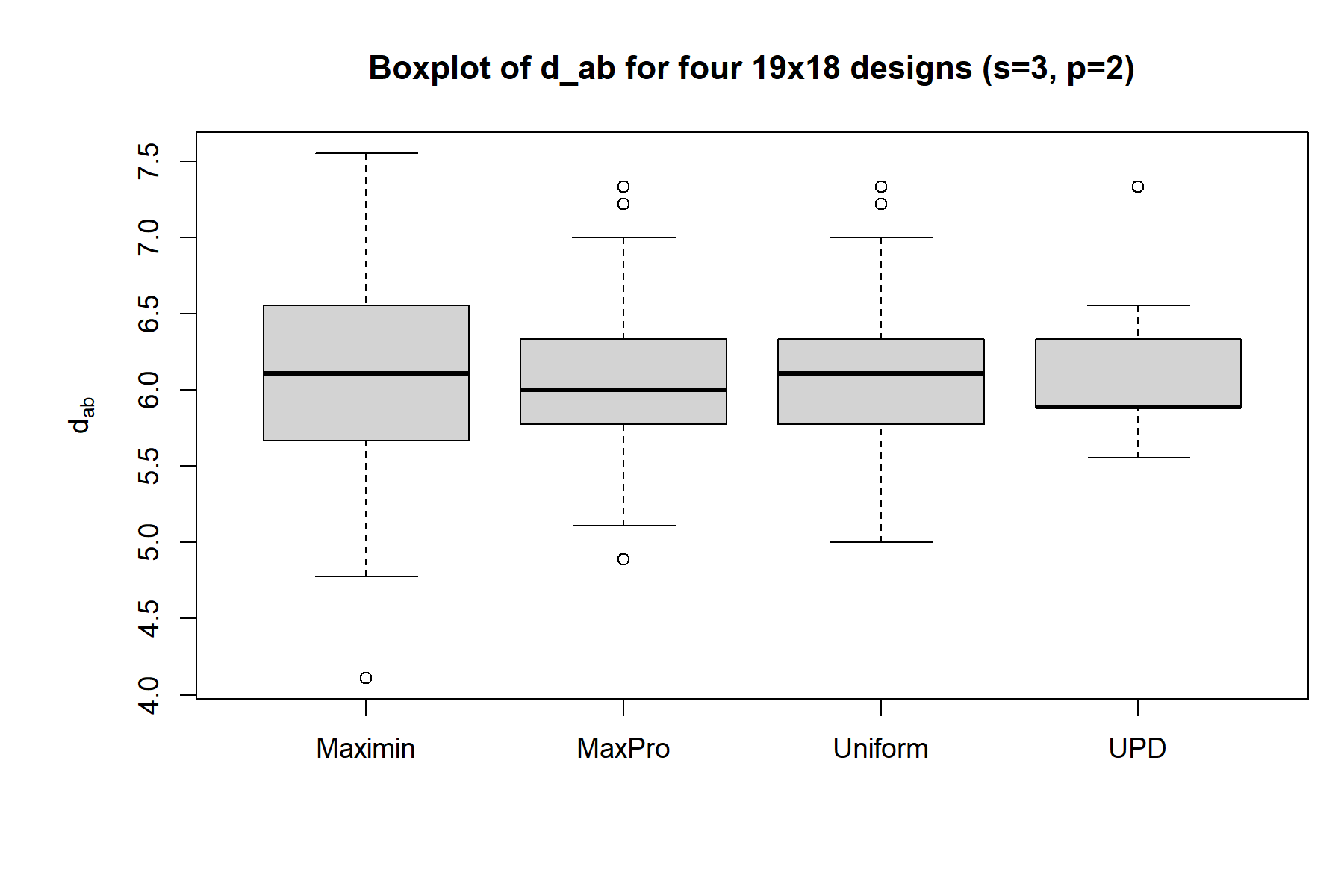}\\
\small (b) $s=3$, $p=2$
\end{minipage}
\caption{Box plots of pairwise weighted hierarchical distances $d_{ab}$ for $19 \times 18$ LHDs.}
\label{fig:sec4_19x18_dab_box}
\end{figure}

Motivated by the row-pairwise distance representation in Theorem~\ref{thm:bounds_PhiSD}, we further compare these designs using the empirical distributions of the pairwise weighted hierarchical distances $d_{ab}$. As shown in Figure~\ref{fig:sec4_19x18_dab_box}, UPD exhibits a more concentrated distance distribution under both kernel settings. This observation is consistent with its favorable $\Phi_{\rm SD}$ values and supports the interpretation of the criterion as a measure of hierarchical space-filling robustness.
Overall, these examples show that $\Phi_{\rm SD}$ is informative even when several full-dimensional criteria provide similar rankings.

The numerical comparisons above use $y=1$ in the exponential weighting scheme, which corresponds to constant weights.
This choice provides a relatively stringent assessment of stratification, since it retains the contribution of finer levels more strongly than choices with $y<1$.
In applications, the value of $y$ should reflect the resolution at which two-factor projection uniformity is to be assessed.
Smaller values of $y$ emphasize coarse stratification in two-factor projections, whereas larger values give more weight to finer partitions.
Thus, a relatively small value of $y$ is suitable when broad space-filling coverage is the main concern, while a larger value, including the constant-weight case $y=1$, may be useful when finer physical resolution or local variation of the response surface is important.
For high-dimensional designs, one may also use the stabilization choice $y=2/(m+1)$ recommended by \citet{TianXu25}.
When there is no clear prior preference, one may evaluate $\Phi_{\rm SD}$ over a small grid of $y\in(0,1]$ and examine whether the ranking of candidate designs is stable.

\section{Further results}\label{sec5}
This section investigates the relationship between the uniform projection criterion $\Phi_{\rm SD}$ and some existing criteria or designs under the space-filling hierarchy framework.

\subsection{An NRT representation of \texorpdfstring{$\Phi_{\rm SD}$}{PhiSD}}


Let $\bo=\{\omega(0),\ldots,\omega(p)\}$ be the weights in the stratified kernel $K^R$ defined in~\eqref{KR}.
Recall the constants $A_0,A_1,B,C$ and $C_{\rm SD}(m,s,p)$ in Theorem~\ref{thmSD2}, and
$A_0^{(\ell)}=\sum_{i=0}^{\ell}\omega(i)/s^i$ for $\ell=0,\ldots,p-1$ in Theorem~\ref{thm:bounds_PhiSD}.
For notational convenience in this section, we also set
\[
A_0^{(p)}:=A_0=\sum_{i=0}^{p}\frac{\omega(i)}{s^i}.
\]
Following \citet{TianXu25}, define for $d=0,\ldots,p$ the partial-sum function
\begin{equation*}
g(d;\bo):=\sum_{i=0}^{p-d}\frac{\omega(i)}{s^i}=A_0^{(p-d)}.
\end{equation*}

Let $\rho(\cdot,\cdot)$ be the NRT-distance on $\mathbb{Z}_{s^p}$ defined in~\eqref{eq:NRT-dist}.
For $x,y\in\mathbb{Z}_{s^p}$, write $\tilde x=(2x+1)/(2s^p)$ and $\tilde y=(2y+1)/(2s^p)$.
Then the key identity used in \citet{TianXu25} can be written in our notation as
\begin{equation}\label{eq:g_delta_identity_sec5}
g\bigl(\rho(x,y);\bo\bigr)
=\sum_{i=0}^p\frac{\omega(i)}{s^i}\,\delta_i(\tilde x,\tilde y)
= A_0^{(p-\rho(x,y))}.
\end{equation}

By Definition~\ref{def:habdist} and \eqref{eq:g_delta_identity_sec5}, for any two runs $a,b\in\{1,\ldots,n\}$,
\begin{equation}\label{eq:dab_as_A0ell_sec5}
\dist_{ab}
=\sum_{k=1}^m\sum_{i=0}^p \frac{\omega(i)}{s^i}\Bigl\{1-\delta_i(z_{ak},z_{bk})\Bigr\}
=\sum_{k=1}^m\Bigl\{A_0-A_0^{(p-\rho(x_{ak},x_{bk}))}\Bigr\}.
\end{equation}
Thus, $\dist_{ab}$ aggregates one-dimensional penalties that depend only on the NRT distance
between $x_{ak}$ and $x_{bk}$.

 \begin{coro}\label{cor:PhiSD_NRT_sec5}
Let $D=(x_{ik})$ be a U-type $(n,m,s^p)$ design. Then we have
\begin{equation}\label{eq:PhiSD_NRT_sec5}
\Phi_{\rm SD}(D)
=\frac{1}{n^2m(m-1)}\sum_{a=1}^n\sum_{b=1}^n
\left(\sum_{k=1}^m\Bigl\{A_0-A_0^{(p-\rho(x_{ak},x_{bk}))}\Bigr\}\right)^2
+ C_{\rm SD}(m,s,p),
\end{equation}
where $C_{\rm SD}(m,s,p)$ is the constant defined in Theorem~\ref{thmSD2}.
\end{coro}

Corollary~\ref{cor:PhiSD_NRT_sec5} makes the hierarchical nature of $\Phi_{\rm SD}$ explicit: larger NRT distance $\rho(x_{ak},x_{bk})$
reduces the partial sum $A_0^{(p-\rho(x_{ak},x_{bk}))}$ and hence increases the penalty
$A_0-A_0^{(p-\rho(x_{ak},x_{bk}))}$ in \eqref{eq:dab_as_A0ell_sec5}.

\subsection{Stratification pattern enumerators and two-dimensional projections}


Following \citet{TianXu22}, for $u,x\in \mathbb{Z}_{s^p}$ define
$\langle u,x\rangle=\sum_{i=1}^p f_{p-i+1}(u)\, f_i(x)$ and let $\xi=\exp(2\pi \sqrt{-1}/s)$.
Define the character $\chi_u(x)=\xi^{\langle u,x\rangle}$.
For $\bu=(u_1,\ldots,u_m)$ and $\bx=(x_1,\ldots,x_m)$ in $\mathbb{Z}_{s^p}^m$, define
\[
\rho(u):=\rho(u,0),\qquad
\rho(\bu)=\sum_{k=1}^m\rho(u_k),\qquad
\chi_{\bu}(\bx)=\prod_{k=1}^m\chi_{u_k}(x_k).
\]
For an $n\times m$ design $D=(x_{ij})$ with levels in $\mathbb{Z}_{s^p}$, let
$\chi_{\bu}(D)=\sum_{\bx\in D}\chi_{\bu}(\bx)$, where the sum is over all rows of $D$.
For $j=0,\ldots,mp$, define
\begin{equation}\label{eq:sp_sec5}
S_j(D)=n^{-2}\sum_{\rho(\bu)=j}\bigl|\chi_{\bu}(D)\bigr|^2.
\end{equation}
The vector $(S_1(D),\ldots,S_{mp}(D))$ is called the \emph{space-filling pattern}.
The minimum aberration-type space-filling criterion ranks designs by sequentially minimizing
$S_j(D)$ for $j=1,\ldots,mp$; moreover, \citet{TianXu22} show that a GSOA has strength $t$
if and only if $S_j(D)=0$ for $1\le j\le t$.

\citet{TianXu24} proposed the \emph{stratification pattern enumerator}
\begin{equation}\label{eq:E_def_sec5}
E(D;y)=\sum_{j=0}^{mp} S_j(D)\,y^j,
\end{equation}
where $y$ is an indeterminate variable.
When $y > 0$ is sufficiently small, $E(D; y)$ preserves the same ranking order as the minimum aberration-type space-filling criterion.
To focus on two-dimensional projection stratification, we generalize $E(D;y)$
by defining the
\emph{two-dimensional projection pattern enumerator}  as
\begin{equation}\label{eq:E2_def_sec5}
E_2(D; y)=\frac{2}{m(m-1)} \sum_{|u|=2}E(D_u; y),
\end{equation}
where $u\subset\{1,2,\ldots,m\}$ with $|u|=2$, and $D_u$ represents the projected design of $D$ onto the factors indexed by $u$.

Set the weights in \eqref{KR} as
$\omega(i)=(s^2y)^i$ for $i=0,\ldots,p-1$ and
$\omega(p)=(s^2y)^p/(1-y)$ with $y\in(0,1)$.
Then \citet[Theorem~4]{TianXu24} established
\begin{equation}\label{eq:SD2_E_sec5}
SD(D)^2=\frac{E(D;y)-1}{(1-y)^m}.
\end{equation}
Applying \eqref{eq:SD2_E_sec5} to each two-factor projection and using the definition of $\Phi_{\rm SD}$ in~\eqref{eqphiDisc}
yields the following identity under the same exponential weighting scheme.

\begin{coro}\label{cor:PhiSD_E2_sec5}
Let $D$ be a U-type $(n,m,s^p)$ design. Under the exponential weighting scheme above, we have
\begin{equation*}
\Phi_{\rm SD}(D)=\frac{E_2(D;y)-1}{(1-y)^2},
\end{equation*}
where $y\in(0,1)$ is an indeterminate variable.
\end{coro}

It follows from Corollary~\ref{cor:PhiSD_E2_sec5} that for two $n\times m$ designs $D_1$ and $D_2$  with levels in $\mathbb{Z}_{s^p},$ if $D_1$ is more space-filling than $D_2$ based on the two-dimensional projection pattern enumerator, then there exists  $\epsilon>0$ such that $\Phi_{\rm SD}(D_1)<\Phi_{\rm SD}(D_2)$ for any $y\in(0,\epsilon),$  provided the uniform projection of the stratified $L_2$-discrepancy uses the exponential weighting scheme above.

\subsection{A remark on SOAs of strength \texorpdfstring{$2+$}{2+}}


The application of $E_2(D; y)$ warrants further exploration. For example, it is closely related to the strong orthogonal arrays (SOAs) of strength $2+$ introduced by \cite{HeChengTang18}. A design $D$ is referred to as a strong orthogonal array of strength $2+$, denoted by $\rSOA(n, m, s^2, 2+)$, if any two columns of $D$ can be collapsed into both an $\rOA(n, 2, s^{2} \times s, 2)$ and an $\rOA(n, 2, s \times s^{2}, 2)$. Such SOAs of strength $2+$ are not members of GSOAs. While they achieve stratifications on $s^2 \times s$ and $s \times s^2$ grids in any two dimensions (similar to GSOAs of strength 3), they do not necessarily achieve stratifications on $s \times s \times s$ grids in any three dimensions.

Let $\bar{S}_1(D)$, $\bar{S}_2(D)$, and $\bar{S}_3(D)$ denote the coefficients of $y$, $y^2$, and $y^3$, respectively, in $E_2(D; y)$, i.e.,
\[
E_2(D;y)=1+\bar S_1(D)y+\bar S_2(D)y^2+\bar S_3(D)y^3+\cdots.
\]
As a two-dimensional analogue of the space-filling pattern, it follows directly that $D$ is an $\rSOA(n, m, s^2, 2+)$ if and only if $\bar{S}_1(D) = \bar{S}_2(D) = \bar{S}_3(D) = 0$.
Thus $(\bar S_1,\bar S_2,\bar S_3,\ldots)$ is well suited for describing pairwise stratification properties.
A systematic study of $E_2(D;y)$ beyond its connection to $\Phi_{\rm SD}$ is an interesting direction for future work.

\section{Discussion}\label{sec:conclusion}

This paper proposes a uniform projection criterion based on the stratified $L_2$-discrepancy to promote
two-dimensional projection uniformity in high-dimensional space-filling designs.
By averaging the stratified discrepancy over all two-factor projections, the criterion mitigates the curse of dimensionality and is consistent with the effect sparsity principle, since many analyses are driven primarily by low-order effects.
We derive an explicit representation of $\Phi_{\rm SD}(D)$ in terms of the weighted hierarchical distances $\dist_{ab}$,
together with sharp lower and upper bounds. These results lead to an $O(n^2m)$ evaluation of $\Phi_{\rm SD}(D)$.
Importantly, we show that designs that are optimal for the full-dimensional stratified $L_2$-discrepancy are also optimal
for the proposed uniform projection criterion, establishing a coherent hierarchy that links low-dimensional projection
uniformity to high-dimensional stratified structure.

A three-dimensional projection analogue can be defined in a similar way, but its general form involves additional interaction terms and is therefore much harder to analyze and optimize. Developing a full theory for three-dimensional stratified projection uniformity, beyond the restricted strength-two OA case discussed after Theorem~\ref{thm:bounds_PhiSD} and Example~\ref{ex2}, is a potential direction for future work.

Our study focuses on symmetric designs, where all factors share the same number of levels. This assumption simplifies the theory, but many applications involve mixed-level (asymmetric) factors.
Extending the stratified projection criterion and the associated theory to asymmetric designs is a meaningful direction for future work.

Finally, \citet{Zhang25} generalized the space-filling pattern of \citet{TianXu22} and introduced a two-dimensional projection pattern
to characterize two-factor stratification.
Their criterion is projection-specific (fixing a given pair of columns) rather than projection-averaged. Therefore, as noted at the end of Section~\ref{sec5}, a comprehensive study of the uniform projection criterion and its relationship with the space-filling patterns (or the stratification pattern enumerator) would be valuable.

\backmatter

 \par

\section*{Supplementary Materials}
All  proofs and additional design tables are provided in the Supplementary Materials.


\end{document}